\documentclass[12pt,twoside,a4paper]{amsart}
\usepackage{amssymb}
\usepackage{graphicx}

\date{\today}

\def\dbar{\bar\partial}

\def\C{{\mathbb C}}
\def\P{{\mathbb P}}

\def\Re{{\rm Re\,  }}

\def\be{\begin{equation}}
\def\ee{\end{equation}}

\def\codim{\text{codim}\,}
\def\ann{\text{Ann}\,}
\def\supp{\text{supp}\,}

\def\scalar{\cdot}
\def\1{\mathbf 1}
\def\0{mathbf 0}
\def\sgn{\text{sgn }\,}

\newtheorem{thm}{Theorem}[section]
\newtheorem{lma}[thm]{Lemma}
\newtheorem{cor}[thm]{Corollary}

\newtheorem{claim}{Claim}

\theoremstyle{definition}

\theoremstyle{remark}

\newtheorem{preremark}{Remark}
\newtheorem{preex}{Example}

\newenvironment{remark}{\begin{preremark}}{\qed\end{preremark}}
\newenvironment{ex}{\begin{preex}}{\qed\end{preex}}

\numberwithin{equation}{section}

\begin{document}

\title[Residue currents of monomial ideals
]{Residue currents of monomial ideals}

\date{\today}

\author{Elizabeth Wulcan}

\address{Department of Mathematics\\Chalmers University of Technology and the University of G\"oteborg\\S-412 96 G\"OTEBORG\\SWEDEN}

\email{wulcan@math.chalmers.se}



\begin{abstract}
We compute residue currents of Bochner-Martinelli type associated with a monomial ideal $I$, by methods involving certain toric varieties. 
In case the variety of $I$ is the origin, we give a complete description of the annihilator of the currents in terms of the associated Newton diagram. In particular, we show that the annihilator is strictly included in $I$, unless $I$ is defined by a complete intersection. We also provide partial results for general monomial ideals.
\end{abstract}


\maketitle

\section{Introduction}
Let $f$ be a tuple of holomorphic functions $f_1,\ldots,f_m$ in $\mathbb C^n$ and let $Y=\{f_1=\ldots =f_m=0\}$. If $f$ is a complete intersection, that is, the codimension of $Y$ is $m$, the duality theorem, due to Dickenstein-Sessa, ~\cite{DS}, and Passare, ~\cite{P}, asserts that a holomorphic function $h$ locally belongs to the ideal $(f)=(f_1,\ldots ,f_m)$ if and only if $h R^f_{CH}=0$, where $R^f_{CH}$ is the Coleff-Herrera residue current of $f$. In ~\cite{PTY}, Passare, Tsikh and Yger introduced residue currents for arbitrary $f$ by means of the Bochner-Martinelli kernel. For each ordered index set $\mathcal I\subseteq \{1,\ldots, m\}$ of cardinality $k$, let $R^f_{\mathcal I}$ be the analytic continuation to $\lambda=0$ of 
\begin{equation*}
\dbar |f|^{2\lambda}\wedge
\sum_{\ell=1}^k(-1)^{k-1}
\frac
{\overline{f_{i_\ell}}\bigwedge_{q\neq \ell} \overline{df_{i_q}}}
{|f|^{2k}},
\end{equation*}
where $|f|^2=|f_1|^2+\ldots+|f_m|^2$.
Then $R^f_\mathcal I$ is a well-defined $(0,k)$-current with support on $Y$, that vanishes whenever $k<\codim Y$ or $k>\min(m,n)$. In case $f$ defines a complete intersection, the only nonvanishing current, $R^f_{\{1,\ldots,m\}}$, is shown to coincide with the Coleff-Herrera current. 

The concept of Bochner-Martinelli residue currents was further developed by Andersson in ~\cite{A}. From his construction, based on the Koszul complex, follows that $hR^f_{\mathcal I}=0$ for all $\mathcal I$ implies that the holomorphic function $h$ belongs to the ideal $(f)$ locally. Thus, letting $\ann R^f$ denote the annihilator ideal, $\{ h \text{ holomorphic}, hR^f_{\mathcal I}=0, \forall \mathcal I\}$, we have that
\begin{equation}\label{flisa}
\ann R^f\subseteq (f).
\end{equation}
The inclusion is strict in general, and thus the currents $R^f_\mathcal I$ do not fully characterize $(f)$ as in the complete intersection case.
Still the ideal $\ann R^f$ is big enough to catch in some sense the ``size'' of $(f)$. Recall that a holomorhic function $h$ belongs locally to the \emph{integral closure} of $(f)$, denoted by $\overline{(f)}$, if $|h|\leq C|f|$ for some constant ~$C$, or equivalently if $h$ fulfills a monic equation $h^r+g_1 h^{r-1}+\ldots + g_r=0$ with $g_i\in (f)^i$ for $1\leq j\leq r$. In ~\cite{PTY} it was proven that $h R^f_\mathcal I=0$ for any $h$ that is locally in the integral closure of $(f)^k$, where $k=|\mathcal I|$, and thus we get
\begin{equation}\label{spisa}
\overline{(f)^\mu}\subseteq \ann R^f,
\end{equation}
where $\mu=\min(m,n)$.
Now, combining 
 ~\eqref{flisa} and ~\eqref{spisa} yields a proof of the Brian\c{c}on-Skoda theorem ~\cite{BS}: 
$\overline{(f)^\mu}\subseteq (f)$.
This motivates us to study the ideal $\ann R^f$. 

In this paper we compute the Bochner-Martinelli currents $R^f_\mathcal I$ in case the generators $f_i$ are monomials and $Y=\{0\}$. Our main result, Theorem ~\ref{main_thm}, gives a complete description of $\ann R^f$ in terms of the Newton diagram associated with the generators. In particular it turns out that $\ann R^f$ depends only on $(f)$, not on the particular choice of generators. 
Also, it follows that we have equality in ~\eqref{flisa} if and only if $(f)$ is a complete intersection and moreover that the inclusion ~\eqref{spisa} is always strict. The proof of Theorem ~\eqref{main_thm}, given in Section ~\ref{leproof}, amounts to computing residue currents in a certain toric variety constructed from the generators, using ideas originally from Varchenko, ~\cite{V}, and Khovanskii, ~\cite{K}. In Section ~\ref{relative} we provide partial results for the case of general monomial ideals. 

\section{Preliminaries and notation}\label{preliminarier}
Let $A$ be a set in $\mathbb Z_+^n$ and let $z^A$ denote the tuple of monomials $\{z^a\}_{a\in A}$, where $z^a=z_1^{a_1}\cdots z_n^{a_n}$ if $a=(a_1,\ldots,a_n)$. The ideal $(z^A)$ admits a nice geometric interpretation as the set $\cup_{a\in A} (a+\mathbb R_+^n)\subset \mathbb R^n$. Indeed, a holomorphic function is in the ideal precisely when its support ($\supp \sum\varphi_a z^a=\{a\in \mathbb Z_+^n, \varphi_a\neq 0\}$) is in 
$\cup_{a\in A} (a+\mathbb R_+^n)$.
The \emph{Newton polyhedron} $\Gamma^+(A)$ of $A$ is defined as the convex hull of $\cup_{a\in A} (a+\mathbb R_+^n)\subset \mathbb R^n$ and the \emph{Newton diagram} $\Gamma(A)$ of $A$ is the union of all compact faces of the Newton polyhedron. For further reference we remark that the set of vertices of the Newton polyhedron is a subset of $A$, see for example ~\cite{Z}.

We will work in the framework from ~\cite{A} and use the fact that the currents $R_{\mathcal I}^f$ appear as the coefficients of full Bochner-Martinelli current introduced there.
We identify $z^A$ with a section of the dual bundle $E^*$ of a trivial vector bundle $E$ over $\mathbb C^n$ of rank $m=|A|$, endowed with the trivial metric. 
If $\{e_a\}_{a\in A}$ is a global holomorphic frame for $E$ and $\{e^*_a\}_{a\in A}$
is the dual frame, we can write $z^A=\sum_{a\in A} z^a e_a^*$. We let $s$ be the dual section $\sum_{a\in A} \bar z^a e_a$ of $z^A$. Also, we fix an ordering of $A$.

Next, we let 
\begin{equation*}
u =
\sum_\ell\frac{s\wedge(\dbar s)^{\ell-1}}{|z^A|^{2\ell}},
\end{equation*}
where $|z^A|^2=\sum_{a\in A}|z^a|^2$,
be the full Bochner-Martinelli form, 
introduced in ~\cite{A2} in order to construct integral formulas with weight factors in a convenient way.
Then $u$ is a smooth section of $\Lambda(E\oplus T_{0,1}^*(\mathbb C^n))$ (where $e_a\wedge d\bar z_i=-d\bar z_i\wedge e_a$), that is clearly well defined outside $Y=f^{-1}(0)$, and moreover 
\begin{equation*}
\dbar|z^A|^{2\lambda}\wedge u
\end{equation*}
has an analytic continuation as a current to $\Re \lambda > -\epsilon$.
The \emph{(full) Bochner-Martinelli residue current} $R^{z^A}$ is defined as the value at $\lambda=0$. Then $R^{z^A}$ has support on $Y$ and $R^{z^A}=R_p+\ldots +R_\mu$, where $p=\codim Y$ and $\mu=\min(m,n)$, and where $R_k\in \mathcal D'_{0,k}(\mathbb C^n,\Lambda^k E)$, by analogy with the fact that the current $R^f_{\mathcal I}$ vanishes if $|\mathcal I|$ is smaller than $p$ or greater than $\mu$. We should remark that Andersson's construction of residue currents, using kernels of  Cauchy-Fantappi\`e-Leray type, works for sections of any holomorphic vector bundle equipped with some Hermitian metric. Observe that in our case (trivial bundle and trivial metric), though, the coefficients of $R^{z^A}$ are just currents of the type $R^f_{\mathcal I}$. Indeed, letting $s_B$ be the section $\sum_{a\in B} \bar z^a e_a$, we can write
$u$ as a sum, taken over subsets $B$ of $A$, of terms
\begin{equation*}
u_B=
\frac{s_B\wedge (\dbar s_B)^{k-1}}{|z^A|^{2n}},
\end{equation*}
where $k$ is the cardinality of $B$.
The corresponding current,
\begin{equation}\label{arbe}
\dbar|z^A|^{2\lambda}\wedge u_B
\end{equation}
evaluated at $\lambda=0$, denoted by $R^{z^A}_B$ or $R_B$ for short, is then merely the current $R^f_{\mathcal I}$ with $\mathcal I$ corresponding to the subset $B$, times the basis element $e_B=\bigwedge_{a\in B} e_a$, where the wedge product is taken with respect to the ordering. Henceforth we will deal with the 
Bochner-Martinelli currents rather then currents $R^f_{\mathcal I}$.

Let us make an observation that will be of further use. If the section $s$ can be written as $\mu s'$ for some smooth function $\mu$ we have the following homogeneity:
\begin{equation}\label{homogen}
s\wedge (\dbar s)^{k-1}=\mu^k s'\wedge (\dbar s)^{k-1},
\end{equation}
that holds since $s$ is of odd degree.

We will use the notation $\dbar [1/f]$ for the value at $\lambda=0$ of $\dbar |f|^{2\lambda}/f$ and analogously by $[1/f]$ we will mean $|f|^{2\lambda}/f|_{\lambda=0}$, that is just the principal value of $1/f$. By iterated integration by parts we have that
\begin{equation}\label{salt}
\int_z
\dbar\Big[\frac{1}{z^p}\Big] \wedge \varphi dz = 
\frac{2\pi i}{(p-1)!}\frac{\partial^{p-1}}{\partial z^{p-1}}\varphi(0).
\end{equation}
In particular, the annihilator of $\dbar [1/z^p]$ is $(z^p)$. The currents $R^{z^A}_B$ will typically be tensor products of currents of this type.

\section{Main results}\label{snuva}
Our main result is an explicit computation of the Bochner-Martinelli residue current $R^{z^A}$ in case $Y$ is the origin. Before stating it let us introduce some notation. We say that a subset $B=\{a_1,\ldots, a_n\}\subseteq A$ is \emph{essential} if there exists a facet $F$ of $\Gamma^+(A)$ such that $B$ lies in $F$ and if in addition $B$ spans $\mathbb R^n$, that is $\det (a_{1},\ldots,a_{n})\neq 0$. It follows, when $Y=\{0\}$, that the essential sets are contained in the Newton diagram $\Gamma(A)$. Indeed, $Y=\{0\}$ precisely when $A$ intersects all axes in $\mathbb Z^n$ and thus the only non-compact faces of $\Gamma^+$ are contained in the coordinate planes in $\mathbb Z^n$. But if $B$ is contained in a coordinate plane, $B$ cannot span $\mathbb R^n$. Also, when $Y=\{0\}$, all points in $A\cap\Gamma(A)$ are in fact contained in some essential set. Next, if $B$ is a subset of $A$, let $\alpha^B=\sum_{a\in B}a$. Notice that if $B$ is essential, then $\alpha^B$ lies on $n\Gamma$. In fact, $\alpha^B/n$ is the barycenter of the simplex spanned by $B$. 
We are now ready to formulate our main theorem.
\begin{thm}\label{main_thm}
Let $z^A, A \subseteq \mathbb Z^n_+$ be a tuple of monomials in $\mathbb C^n$ such that $\{z^A=0\}=\{0\}$, and let $R^{z^A}$ be the corresponding Bochner-Martinelli residue current . Then 
\begin{equation*}
R^{z^A}=\sum_{B\subseteq A} R_B,
\end{equation*}
where 
\begin{equation}\label{muffel}
R_B= 
C_B ~~\dbar \Big [\frac{1}{z_1^{\alpha^B_1}} \Big ]\wedge\ldots\wedge
 \dbar \Big[\frac{1}{z_n^{\alpha^B_n}}\Big]\wedge e_B,
\end{equation}
and where $C_B$ is a constant that is nonzero if $B$ is an essential set and zero otherwise.
\end{thm}
An immediate consequence is that if $B$ is essential then
\begin{equation*}
\ann R_B=(z_1^{\alpha^B_1},\ldots,z_n^{\alpha^B_n}),
\end{equation*}
where $\ann R_B$ just denotes the ideal of holomorphic functions annihilating $R_B$.
Note in particular that $\ann R_B$ depends only on the set $B$ and not on the remaining $A$. Furthermore, since the basis elements
$e_B$ are all different it follows that 
\[
\ann R^{z^A}=\bigcap_{B \text{ essential}} \ann R_B.
\]
Thus, $\ann R^{z^A}$ is fully determined by the Newton diagram $\Gamma(A)$ and the points in $A$ lying on it. In particular $\ann R^{z^A}$ depends only on the ideal, not on the particular choice of generators. We also see that different monomial ideals $(z^A)$ and $(z^{A'})$ give rise to the same annihilator ideal if and only if $A\cap \Gamma(A)=A'\cap \Gamma(A')$.

Furthermore, Theorem ~\ref{main_thm} implies that the inclusion ~\eqref{flisa} is strict unless we have a complete intersection.
\begin{thm}\label{main_corollary}
Let $z^A, A\subseteq \mathbb Z^n_+$, be a tuple of monomials such that $\{z^A=0\}=\{0\}$, and let 
$ R^{z^A}$ be the corresponding Bochner-Martinelli residue current. Then
\begin{equation}\label{orm}
\ann R^{z^A}=(z^A) 
\end{equation}
if and only if $(z^A)$ can be generated by a complete intersection.
\end{thm}
For the proof we need a simple lemma.
\begin{lma}\label{kaffi}
Let $B$ be an essential subset of $A$ such that $(z^B)\subseteq \ann R_B$. 
Then $(z^B)$ is a complete intersection.
\end{lma}
\begin{proof}
Denote the elements in $B$ by $a_i, i=1,\ldots,n$ and let $\geq$ be the natural partial order on $\mathbb Z^n$.
Suppose that $(z^B)\subseteq \ann R_B$. We have that $z^{a_i}\in \ann R_B$ precisely when one of the generators of $\ann R_B$ divides $z^{a_i}$, that is, when
\begin{eqnarray*}
(a_{1i},\ldots,a_{ni})&\geq & (\sum_j a_{1j},0,\ldots,0) \text{ or }\\
(a_{1i},\ldots,a_{ni})&\geq & (0,\sum_j a_{2j},0,\ldots,0) \text{ or }\\
&\vdots&\\
(a_{1i},\ldots,a_{ni})&\geq & (0,\ldots, 0,\sum_j a_{nj}).
\end{eqnarray*}
This set of inequalities holds for all $1\leq i\leq n$, and it is easy to see that this implies first that $a_{k\ell}\neq 0$ for at most one $k$, which means that $a_\ell$ lies in one of the coordinate axes, and second that there is at least one $a_\ell$ intersecting each coordinate axis. Thus, $B$ intersects all coordinate axes in $\mathbb Z^n$, which in turn implies that $(z^B)$ is a complete intersection.
\end{proof}

\begin{proof}[Proof of Theorem ~\ref{main_corollary}]
We need to show the ``only if'' direction. Suppose that $(z^A)= \ann R^{z^A}$ and let $B$ be an essential subset. Clearly essential subsets always exist, since otherwise $R^{z^A}=0$ and $\ann R^{z^A}=(z^A)$ is the whole ring of holomorphic functions, which contradicts that $Y=\{0\}$. Now, in particular $(z^B)\subseteq \ann R_B$, and by Lemma ~\ref{kaffi}, $(z^B)$ is a complete intersection. Thus  
\[(z^B)=\ann R^{z^B}=\ann R^{z^A}_B\supseteq \ann R^{z^A}=(z^A)\supseteq(z^B),\]
where the second equality follows since $\ann R_B$ only depends on $B$ and not on $A$. Hence $(z^A)=(z^B)$ and the result follows.
\end{proof}

We give some examples to illustrate Theorems ~\ref{main_thm} and ~\ref{main_corollary}.
\begin{ex}\label{basex}
Let 
\[A=\{ a^1=(8,0), a^2=(6,1), a^3=(2,3), a^4=(1,5), a^5=(0,6)\}\subseteq \mathbb Z^2.\]
We identify the ideal $(z^A)$ with the set $\bigcup_{a\in A}(a+\mathbb R_+^n)$ as in Figure 1, where we have also depicted the Newton diagram $\Gamma$. Such pictures of monomial ideals are usually referred to as \emph{staircase diagrams}, see ~\cite{MS}. The points in $A$ should be recognized as the ``inner corners'' of the staircase. 
\begin{figure}\label{staircase}
\begin{center}
\includegraphics{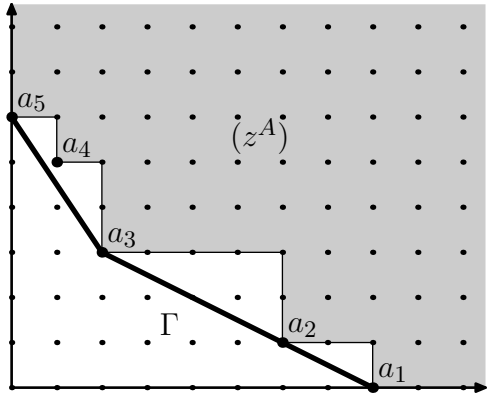}
\caption{The ideal $(z^A)$ and the Newton diagram $\Gamma(A)$ in Example ~\ref{basex}}
\end{center}
\end{figure}
The Newton diagram $\Gamma(A)$ consists of two facets, one with vertices $a^1$ and $a^3$ and the other one with vertices  $a^3$ and $a^5$, and thus we have the essential sets
\[
\{a^1,a^2\},\{a^1,a^3\},\{a^2,a^3\},\{a^3,a^5\},
\]
with 
\[
\alpha^{12}=(14,1), \alpha^{13}=(10,3),\alpha^{23}=(8,4),\alpha^{35}=(2,9),
\]
respectively.
It follows from Theorem ~\ref{main_thm} that 
\[
\ann R^{z^A}=
(z_1^{14},z_2)\cap (z_1^{10},z_2^3)\cap (z_1^8,z_2^4)\cap (z_1^2,z_2^9),
\]
which is equal to the ideal $(z_1^{14},z_1^{10}z_2,z_1^8z_2^3,z_1^2z_2^4,z_2^9)$,
see Figure ~2.
Observe that $\ann R^{z^A}$ is given by the staircase diagram with $\alpha^{ij}$ as ``outer corners''. Note also that $\ann R^{z^A}$ does not depend on $a^4$, which lies in the interior of $\Gamma^+(A)$.
\begin{figure}\label{banni}
\begin{center}
\includegraphics{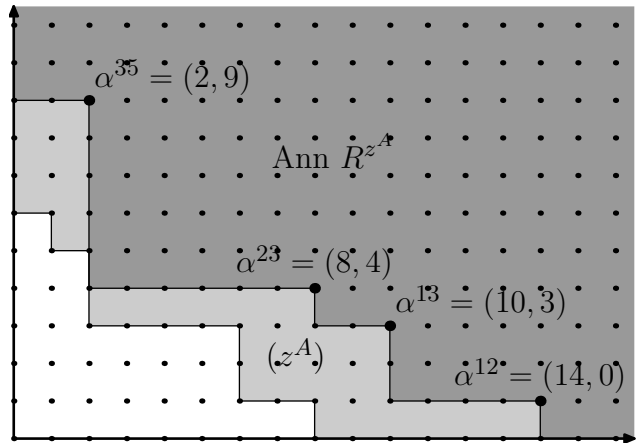}
\caption{The ideals $\ann R^{z^A}$ (dark grey) and $(z^A)$ (light grey) in Example ~\ref{basex}}
\end{center}
\end{figure}
\end{ex}
\begin{ex}\label{fullstex}
Consider the complete intersection $\{z_1^{a^1},\ldots,z_n^{a^n}\}$. The associated Newton diagram is the $n$-simplex spanned by 
\[
A=\{(a^1,0,\ldots,0),\ldots,(0,\ldots,0,a^n)\}
\]
and there exists only one essential set, namely $A$ itself, with $\alpha^A=(a^1,\ldots,a^n)$. Thus according to Theorem ~\ref{main_thm},
\[
\ann R^{z^A}=(z_1^{a^1},\ldots,z_n^{a^n}),
\]
so the annihilator ideal is equal to $(z^A)$, which we already knew. Figure ~3 illustrates the two ways of thinking of the ideal when $n=2$; either as a staircase with $(a^1,0)$ and $(0,a^2)$ as inner corners or as a staircase with $\alpha^A=(a^1,a^2)$ as the (only) outer corner.
\begin{figure}\label{kruka}
\begin{center}
\includegraphics{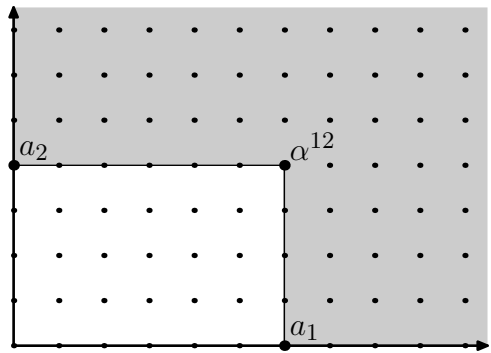}
\caption{A complete intersection}
\end{center}
\end{figure}
\end{ex}
\begin{ex}\label{palm}
We should remark that not all monomial ideals arise as annihilator ideals associated with monomial ideals. The idea is that the outer corners of the staircase of an annihilator ideal must lie on a hypothetical Newton diagram. Indeed, from the discussion just before Theorem ~\ref{main_thm} we know that each $\alpha^B$ corresponding to an essential set $B$ lies on $n\Gamma$. In other words, the lines joining adjacent outer corners must lie on the boundary of a convex domain above the staircase, and thus a necessary condition is that the ``slope'' of the staircase decreases while we are descending it.

For example, consider the ideal
\[
I=(z_1^5, z_1^4z_2^2, z_1z_2^4, z_2^5)
\]
with staircase diagram as in Figure 4, where we have also marked the slope.
\begin{figure}\label{psalm}
\begin{center}
\includegraphics{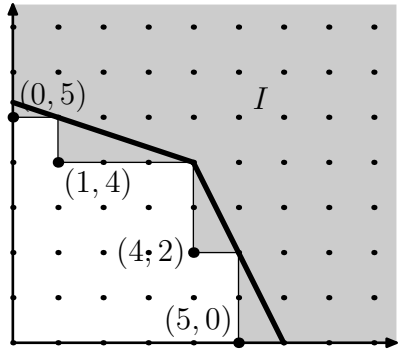}
\caption{The ideal in Example ~\ref{palm}. The thick lines illustrate the ``slope'' of the staircase.}
\end{center}
\end{figure}
Clearly, the outer corners cannot lie on the boundary of a convex Newton polyhedron, and thus $I$ is not an annihilator ideal.
\end{ex}

\begin{remark}
Observe that adding an extra generator to an ideal $(z^A)$ does not necessarily make the corresponding annihilator ideal smaller or larger. However, with a fixed Newton diagram an extra generator can only make the annihilator ideal smaller. In fact, given $\Gamma$, $\ann R^{z^A}$ is maximal if $A$ is chosen as the vertex set of $\Gamma$ and minimal if $A$ is all integer points on $\Gamma$, as we will see in Example ~\ref{hela}.
\end{remark}

Let us now consider the inclusion ~\eqref{spisa}. We start by interpreting the left hand side in case $f$ is monomial. First, we make the following observation.
\begin{lma}\label{loflof}
The integral closure of the monomial ideal $(z^A)$ is the monomial ideal generated by $z^a, a\in \Gamma^+(A)$. 
\end{lma}
The result is well known from algebraic contexts, see for example ~\cite{T}. We supply a proof, however, using the analytic definition of integral closure.
\begin{proof}
We start by proving that $z^b\in\overline{(z^A)}$ for any $b\in \mathbb Z^n\cap \Gamma^+(A)$. It suffices to show that 
\begin{equation}\label{mjolk}
|z^b|\leq C |z^A|=C\sum_{a\in A}|z^a|,
\end{equation}
if $b\in\Gamma(A)$. Indeed, if $b\in \Gamma^+(A)$, we have that $b=c~b'$ for some $b'\in \Gamma(A)$ and $0<c\leq 1$ and so $|z^b|\leq |z^{b'}|$.
Suppose that $b$ lies on the facet $F$ spanned by $a^1,\ldots,a^n$. Then $b=\sum_{i=1}^n\lambda_i a^i$ for some $\lambda_i\geq 0$ such that $\sum_{i=1}^n\lambda_i=1$ and thus $|z^b|=|z^{\sum_{i=1}^n\lambda_i a^i}|=\prod |z^{a^i}|^{\lambda_i}$. Observe that $\prod_{i=1}^n x_i^{\lambda_i}\leq \sum x_i$ if $x_i\geq 0$, $\lambda_i\geq 0$ and $\sum_{i=1}^n\lambda_i=1$. To see this, take the logarithm of each side and use that it is a concave function. Thus ~\eqref{mjolk} follows (and we can choose $C$ to be ~$1$).

Conversely, we need to show that ~\eqref{mjolk} cannot hold if $b\notin \Gamma^+(A)$. Notice that such a $b$ equals $cb'$ for some $b'\in \Gamma(A)$ and $c> 1$. Suppose that $b'$ lies on the facet $F$ with non-negative normal direction $\rho$ and observe that $\rho\scalar a\geq \rho\scalar b'$ for all $a\in A$. Now, for $s\in \mathbb R$ choose $z(s)\in \mathbb C^n$, such that $|z_i(s)|=\exp(s\rho_i)$. Then 
\begin{equation*}
\frac{|z(s)^A|}{|z(s)^b|}=
\frac{\sum_{a\in A} \exp(s\rho\scalar a)}{\exp(cs\rho\scalar b')}\leq
|A| \exp(s \rho\scalar b'(c-1))\to 0 
\end{equation*}
when $s\to -\infty$ and thus ~\eqref{mjolk} cannot hold.
\end{proof}
Next, we claim that the ideal $\overline{(z^A)^r}$ is generated by $z^a, a \in r\Gamma^+(A)$. The ideal $(z^A)^r$ is generated by $z^a, a\in A+\ldots +A$ ($r$ times), so we need to show that the Newton polytope of $A+\ldots +A$ is equal to $r\Gamma^+(A)$. But $A+\ldots +A\supseteq rA$ and thus $\Gamma^+(A+\ldots + A)\supseteq \Gamma^+(rA)=r\Gamma^+(A)$. On the other hand $A+\ldots+ A\subseteq \Gamma^+(A)+\ldots+\Gamma^+(A)=r\Gamma^+(A)$, where the equality holds since $\Gamma^+(A)$ is a convex set, and so it follows that $\Gamma^+(A+\ldots +A)\subseteq r\Gamma^+(A)$.
\begin{cor}\label{yes}
Suppose $n\geq 2$. Let $z^A$ be as in Theorem ~\ref{main_thm}.
Then the integral closure of the ideal $(z^A)^n$ is strictly included in $\ann R^{z^A}$.
\end{cor}
Observe that Corollary ~\ref{yes} fails when $n=1$. Then, in fact, $(z^A)=\ann R^{z^A}= \overline{(z^A)}$.
\begin{proof}
Let $(b_1,0,\ldots,0)$ be the intersection between $\Gamma(A)$ and the $x_1$-axis and let $f=z_1^{nb_1-1}$. Then $(nb_1-1,0,\ldots,0)\notin n\Gamma^+(A)$ and thus $f\notin\overline{(z^A)^n}$. However, $f\in \ann R_B$ for all essential $B$. To see this, observe that the simplex spanned by the intersection points between $\Gamma$ and the axes separates $\Gamma$ from $\{x_1=b_1\}$, and so $\Gamma$ intersects the hyperplane $\{x_1=b_1\}$ only at the point $(b_1,0,\ldots,0)$. This implies in particular that $\alpha_1^B\leq nb_1-(n-1)$ for all essential $B$ and thus $f\in (z_1^{\alpha^B_1})\subseteq \ann R_B$. Hence we have found a function $f$ in $\ann R^{z^A}\setminus \overline{(z^A)^n}$.
\end{proof}
Another, probably more illuminating, way of thinking of the ideals is in terms of staircase diagrams as in the examples above. The fact that the ideal $\overline{(z^A)^n}$ is generated by $\{z^a\}, a\in n\Gamma^+$ means that its staircase lies just above $n\Gamma$. On the other hand we know that the outer corners of the staircase of $\ann R^{z^A}$, the $\alpha^B$, lie on $n\Gamma$ and therefore the staircase must lie under $n\Gamma$. Thus the staircase of $\ann R^{z^A}$ is ``strictly lower'' than the staircase of $\overline{(z^A)^n}$ and so the corresponding inclusion of ideals is strict. For an illustration, see Figure ~5, where we have drawn the staircases of the three ideals $(z^A)$, $\ann R^{z^A}$ and $\overline{(z^A)^n}$ for $A$ from Example ~\ref{basex}.

\begin{figure}\label{honung}
\begin{center}
\includegraphics{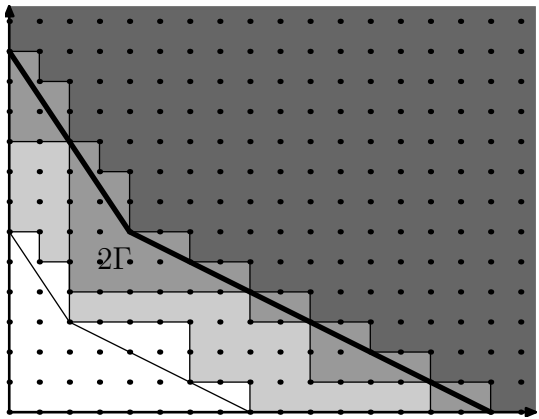}
\caption{The ideals $(z^A)$ (light gray) $\ann R^{z^A}$ (medium gray) and $\overline{(z^A)^2}$ (dark gray) in Example ~\ref{basex}}
\end{center}
\end{figure}

\begin{ex}\label{hela}
Let $\Gamma$ be the simplex with vertices $(3,0)$ and $(0,3)$. In Figure 6 we have drawn the staircases of the ideals $(z^A)$ (light gray), $\ann R^{z^A}$ (medium gray) and $\overline{(z^A)^2}$ (dark gray) for different $A=A_i$ with $\Gamma$ as Newton diagram; more precisely for $A_1=\{(3,0),(0,3)\}$, $A_2=\{(3,0),(2,1),(0,3)\}$, and finally for $A_3=\{(3,0),(2,1),(1,2),(0,3)\}$.
\begin{figure}\label{sko}
\begin{center}
\includegraphics{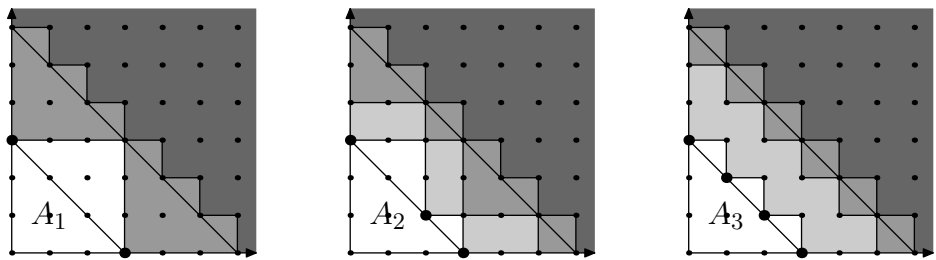}
\caption{The various ideals in Example ~\ref{hela}}
\end{center}
\end{figure}
We see that $\ann R^{z^A}$ decreases when we add points to $A$. In particular integrally closed ideals, that is ideals $I$ such that $\overline I=I$, have the smallest annihilator ideals.
\end{ex}

\section{Proof of Theorem ~\ref{main_thm}}\label{leproof}
The proof of Theorem ~\ref{main_thm} is very much inspired by the the proof of Lemma 2.2 in ~\cite{PTY} and the proof of Theorem 1.1 in ~\cite{A}. We will compute $R^{z^A}$ as the push-forward of a corresponding current on a certain toric variety $\mathcal X$ constructed from the Newton polyhedron $\Gamma^+(A)$. To do this we will have use for the following simple lemma which is proven essentially by integration by parts.
\begin{lma}\label{main_lemma}
Let $v$ be a strictly positive smooth function in $\C$,  
 $\varphi$ a test function in $\C$, and 
$p$  a positive integer.  Then
\[
\lambda\mapsto\int v^\lambda |s|^{2\lambda}\varphi(s)\frac{ds\wedge d\bar s}{s^p}
\]
and
\[
\lambda\mapsto\int \dbar (v^\lambda |s|^{2\lambda})
\wedge \varphi(s)\frac{ds}{s^p}
\]
both have meromorphic  continuations  to the entire plane with poles at rational points on the negative real axis. At $\lambda=0$ they are both independent of $v$ and the second one only depends on the germ of $\varphi$ at the origin. 
Moreover, if $\varphi(s)=\bar s\psi(s)$ or $\varphi=d\bar s\wedge\psi$, then the value of the second integral at $\lambda=0$ is zero.
\end{lma}
Throughout this section we will write $\1$ for the unit vector $(1,1,\ldots,1)$. We will regard the elements in $A$ as column vectors and denote by $B$ the matrix with the vectors in the set $B$ as columns. Also we will use the notation $\widehat{\alpha_i}$ for $\alpha_1\wedge\ldots\wedge\alpha_{i-1}\wedge\alpha_{i+1}\wedge\ldots\wedge\alpha_n$. 

Let us start by describing $\mathcal X$, following ~\cite{BGVY}. Let $\mathcal S$  be the set of normal directions to the facets of $\Gamma^+$ represented by vectors $\rho$ with minimal integer non-negative coefficients. Then $\mathcal S$ provides a partition of the first orthant of $\mathbb R^n$ into a finite number of distinct $n$-dimensional cones. Such a system of cones with the same apex together with their faces is called a \emph{fan}. We say that the fan is generated by $\mathcal S$ and we denote it by $\Delta(\mathcal S)$. By techniques due to Mumford et al., ~\cite{KKMS}, $\mathcal S$ can be completed into a system $\mathcal{\widetilde{S}}$ of vectors $\rho$ such that if $\rho_1,\ldots,\rho_n$ generate one of the $n$-dimensional cones of $\Delta(\mathcal{\widetilde{S}})$, then $\det(\rho_1,\ldots,\rho_n)=\pm 1$. Such a fan is called \emph{regular}.
We will construct $\mathcal X$ by glueing together different copies of $\mathbb C^n$, one for each $n$-dimensional cone of $\Delta(\mathcal{\widetilde{S}})$. Let $\tau$ be such a cone and denote its generators by $\rho_1,\ldots,\rho_n$. Let $\mathcal U$ be the corresponding copy of $\C^n$ with local coordinates $t=(t_1,\ldots,t_n)$. Let $P$ be the matrix with $\rho_i=(\rho_{1i},\ldots,\rho_{ni})$ as rows and let $\Pi$ be the mapping
\begin{eqnarray*}
\Pi:\mathcal U &\to& \mathbb C^n\\
 t & \mapsto & t^P,
\end{eqnarray*}
where $t^P$ is a shorthand notation for
$(t_1^{\rho_{11}}\cdots t_n^{\rho_{n1}}, \ldots, t_1^{\rho_{1n}}\cdots t_n^{\rho_{nn}})$.

Two points  $t\in \mathcal U$ and $t'\in \mathcal U'$ are identified if the monodial map $\Pi'^{-1}\circ \Pi: \mathcal U\to \mathcal U'$ is defined at $t$ and maps $t$ to $t'$. Glueing the charts $\mathcal U$ together induces a proper map $\widetilde\Pi:\mathcal X\to \mathbb C^n$ that is biholomorphic from $\mathcal X \setminus \widetilde\Pi^{-1}(\{z_1\cdots z_n=0\})$ to $\mathbb C^n \setminus \{z_1\cdots z_n=0\}$, that is, outside the coordinate planes. It holds that $\widetilde \Pi^{-1}(\{z_1\cdots z_n=0\})$ is a set of measure zero in $\mathcal X$, and moreover $\widetilde \Pi^{-1}(0)$ consists of a system of various $\C\P^{n-i}$, corresponding to $i$-dimensional cones of the fan $\Delta(\widetilde{\mathcal S})$. In particular, each vector $\rho$, that generates a $1$-dimensional cone, corresponds to a $\C\P^{n-1}$, denoted by $S_\rho$ and obtained by glueing together parts of the charts from the cones determined by $n$-dimensional cones that have $\rho$ as one of its generators. In fact, if the vector $\rho$ determines the coordinate $t_1$ in $\mathcal U$, then $S_\rho$ is covered by the $\{t_1=0\}$-part of $\mathcal U$.

Observe that $R^{z^A}=R_n$ since $Y=\{0\}$. Therefore, we only need to compute the currents $R_B$ when $B$ is a subset of cardinality $n$. For $\Re\lambda$ large enough, ~\eqref{arbe} is integrable and since $\widetilde \Pi$ is biholomorphic outside a set of measure zero it holds that
\begin{equation*}
\int_{\mathbb C^n}
\dbar|z^A|^{2\lambda}\wedge u_B\wedge\phi=
\int_{\mathcal X}
\widetilde\Pi^*(\dbar|z^A|^{2\lambda}\wedge u_B)\wedge\widetilde\Pi^*\phi,
\end{equation*}
if $\phi$ is a test form of bidegree $(n,0)$. It is easy to see that the analytic continuation to $\Re \lambda>-\epsilon$ of $\widetilde\Pi^*(\dbar|z^A|^{2\lambda}\wedge u_B)$ exists in each chart $\mathcal U_\tau$; we will actually compute it below. Thus, because of the uniqueness of analytic continuations,
\begin{equation*}
\widetilde R_B:=\widetilde\Pi^*(\dbar|z^A|^{2\lambda}\wedge u_B)|_{\lambda=0}
\end{equation*}
defines a (globally defined) current on $\mathcal X$ such that $\widetilde \Pi_* \widetilde R_B=R_B$. We will start by computing $\widetilde R_B$ in a fixed chart $\mathcal U_0$ parametrized by $\Pi$ corresponding to the cone $\tau_0$.
\begin{claim}\label{claim1} 
The current 
$\widetilde R_B$ vanishes in $\mathcal U_0$ whenever $B$ is not contained in a facet whose normal direction is one of the generators of ~$\tau_0$. Moreover $\widetilde R_B$ vanishes if $\det B=0$.
\end{claim}
In particular, a necessary condition for $\widetilde R_B$ not to vanish is that $B$ is essential.
\begin{proof}
First, note that the pullback $\Pi^*$ transforms the exponents of monomials by the linear mapping $P$;
\begin{equation}\label{fladder}
\Pi^* z^a=\Pi ^* z_1^{a_1}\cdots z_n^{a_n}=
t_1^{\rho_{11}a_1+\ldots+\rho_{1n}a_n}\cdots
t_n^{\rho_{n1}a_1+\ldots+\rho_{nn}a_n}=
t^{P a}.
\end{equation}
It is well known that for some $a_0\in A, \Pi^* z^{a_0}$ divides $\Pi^* z^a$ for all $a\in A$, and moreover, in view of ~\eqref{fladder} one easily checks that $a_0$ has to be a vertex of $\Gamma^+(A)$. Using this we can write
\[\Pi^* s=\bar t^{Pa_0}s',\]
where $s'$ is the nonvanishing section
\[s'=\sum_{a\in A}\bar t^{P(a-a_0)}e_a,\]
and furthermore
\begin{equation*}
\Pi^*|z^A|^2=|t|^{2Pa_0} \nu(t),
\end{equation*}
where
\[\nu(t)= \sum_{a\in A} |t|^{2P(a-a_0)n}\]
is nonvanishing.
By homogeneity, see ~\eqref{homogen},
\[\Pi^*(s\wedge(\dbar s)^{n-1})=\bar t^{nPa_0} s'\wedge (\dbar s')^{n-1},\]
and thus
\begin{equation}\label{skar}
\widetilde R_B=
\dbar(|t|^{2\lambda Pa_0}\nu^\lambda)~
\frac{s_B'\wedge(\dbar s_B')^{n-1}}
{t^{nPa_0}\nu(t)^n}\bigg |_{\lambda=0}.
\end{equation}
By Leibniz' rule and Lemma ~\ref{main_lemma}, ~\eqref{skar} is equal to a sum of currents
\begin{equation}\label{lila}
\dbar 
\bigg [\frac{1}{t_i^{n\rho_i\scalar a_0}}\bigg ]
\otimes
\bigg [\frac{1}{\prod_{j\neq i}
t_j^{n\rho_j\scalar a_0}}\bigg ]
\wedge 
\frac{s_B'\wedge(\dbar s'_B)^{n-1}}
{\nu(t)^n}.
\end{equation}
We need to compute $s_B'\wedge(\dbar s_B')^{n-1}$. Denote the elements in $B$ by $b_1,\ldots,b_n$ in such a way that
$e_B=e_{b_n}\wedge\ldots\wedge e_{b_1}$. 
Furthermore, let $C$ be the matrix with columns $Pb_i-Pa_0$ so that
\[
s'_B=\sum_i \bar t_1^{c_{1i}}\cdots \bar t_n^{c_{ni}} e_{b_i},
\]
and let $D_i$ be the determinant of $C$ with row $i$ replaced with the unit vector $\1$. 
Then we have the following lemma.
\begin{lma}\label{purple}
We have that
\begin{equation}\label{bubu}
s_B'\wedge (\dbar s_B')^{n-1}
= (n-1)!\bar t^{C\1}\sum_i (-1)^{i-1}~ D_i 
\widehat{\frac{d\bar t_i}{\bar t_i}}
\wedge e_B,
\end{equation}
where
\[
\widehat{\frac{d\bar t_i}{\bar t_i}}=
\frac{d\bar t_1}{\bar t_1}\wedge\ldots\wedge
\frac{d\bar t_{i-1}}{\bar t_{i-1}}\wedge
\frac{d\bar t_{i+1}}{\bar t_{i+1}}\wedge
\ldots\wedge
\frac{d\bar t_{n}}{\bar t_{n}}.
\]
\end{lma}
Observe that all $\bar t_i$ in the denominator are cancelled since ~\eqref{bubu} is in fact smooth.
\begin{proof}
Let $\alpha_j=\bar t_1^{c_{1i}}\cdots \bar t_n^{c_{ni}} e_{b_j}$ and $\beta_i=\frac{d\bar t_i}{\bar t_i}$. 
Then $ s_B'=\sum_{j=1}^n \alpha_j$
and 
\[
\dbar s_B'=\sum_{j=1}^n \sum_{i=1}^n
c_{ij} 
\frac{d\bar t_i}{\bar t_i} \wedge \bar t_1^{c_{1i}}\cdots \bar t_n^{c_{ni}} e_{b_j}=
\sum_{j=1}^n \sum_{i=1}^n c_{ij} \beta_i\wedge \alpha_j.
\]
Thus we get
\begin{multline*}
s_B'\wedge (\dbar s_B')^{n-1}=
\sum_{j=1}^n\alpha_j\wedge(\sum_{j=1}^n\sum_{i=1}^n c_{ij} \beta_i\wedge \alpha_j)^{n-1}=\\
\sum_{\sigma\in S^n}\sum_{\tau\in S^n}
c_{\sigma(2)\tau(2)}\cdots c_{\sigma(n)\tau(n)}
\alpha_{\tau(1)}\wedge\beta_{\sigma(2)}\wedge\alpha_{\tau(2)}\wedge\ldots\wedge
\beta_{\sigma(n)}\wedge\alpha_{\tau(n)}=\\
\sum_{\sigma\in S^n}\sum_{\tau\in S^n}
c_{\sigma(2)\tau(2)}\cdots c_{\sigma(n)\tau(n)}
\beta_{\sigma(2)}\wedge\ldots\wedge\beta_{\sigma(n)}\wedge
\alpha_{\tau(n)}\wedge\ldots\wedge\alpha_{\tau(1)}=\\
\sum_{\sigma\in S^n}\sum_{\tau\in S^n}
(-1)^{\sgn \tau}
c_{\sigma(2)\tau(2)}\cdots c_{\sigma(n)\tau(n)}
\beta_{\sigma(2)}\wedge\ldots\wedge\beta_{\sigma(n)}\wedge
\alpha_{n}\wedge\ldots\wedge\alpha_{1}=\\
\sum_{i=1}^n\sum_{\sigma\in S^n; \sigma(1)=i}
D_i(-1)^{\sgn \sigma}\beta_{\sigma(2)}\wedge\ldots\wedge\beta_{\sigma(n)}\wedge
\alpha_{n}\wedge\ldots\wedge\alpha_{1}=\\
\sum_{i=1}^n
(n-1)!D_i(-1)^{i-1}\beta_1\wedge\ldots \widehat{\beta_i}\ldots \wedge\beta_n\wedge\alpha_{n}\wedge\ldots\wedge\alpha_{1}=\\
(n-1)!\bar t^{C\1}\sum_{i=1}^n (-1)^{i-1}~ D_i 
\widehat{\frac{d\bar t_i}{\bar t_i}}
\wedge e_B.
\end{multline*}
Here $S^n$ just denotes the set of permutations of $\{1,\ldots,n\}$.
\end{proof}
Now 
~\eqref{lila} is equal to
\begin{equation}\label{dennis}
\dbar 
\bigg [\frac{1}{t_i^{n\rho_i\scalar a_0}}\bigg ]
\otimes
\bigg [\frac{1}{\prod_{j\neq i} t_j^{n\rho_j\scalar a_0}}
\bigg ]
\wedge 
\frac{(n-1)!~ D_i ~\bar t^{C\1}} {\nu(t)^n}
\widehat{\frac{d\bar t_i}{\bar t_i}}
\wedge e_B,
\end{equation}
that can vanish for two reasons.
First, by Lemma ~\ref{main_lemma}, ~\eqref{dennis} vanishes whenever the numerator contains a factor $\bar t_i$, that happens
if $c_{ij}>0$ for some ~$j$, which means that
$Pb_j$ has a greater $t_i$-coordinate than $Pa_0$.
Thus, a necessary condition for ~\eqref{dennis} not to vanish is that $P(B)$ is contained in the facet of $P(\Gamma^+)$ parallel to the coordinate plane $\{t_i=0\}$; in other words, since $P$ is invertible, that $B$ is contained in the facet $F_{i}$ of $\Gamma^+$ with normal direction $\rho_i$. Hence the first part of Claim ~\ref{claim1} follows.

Second, ~\eqref{dennis} vanishes if $D_i=0$. 
Assume for simplicity that $i=1$. Then $\rho_1\scalar a$ is constant and equal to $\rho_1\scalar a_0$ on $F_{1}$, that is, $(PB)_{1j}=(Pa_0)_1$ for all $j$, and we get
\begin{displaymath}
D_1=
\begin{array}{|c c c|}
1 & \cdots & 1\\
c_{21} & \cdots & c_{2n}\\
\vdots && \vdots \\
c_{n1} & \cdots & c_{nn}
\end{array}
=
\begin{array}{|c c c|}
 1 & \cdots & 1\\
 (PB)_{21}-(Pa_0)_2 & \cdots & (PB)_{2n}-(Pa_0)_2\\
\vdots && \vdots \\
(PB)_{n1}-(Pa_0)_n & \cdots & (PB)_{nn}-(Pa_0)_n
\end{array}=
\end{displaymath}
\begin{displaymath}
\begin{array}{|c c c|}
1 & \cdots & 1\\
(PB)_{21} & \cdots & (PB)_{2n}\\
\vdots && \vdots \\
(PB)_{n1} & \cdots & (PB)_{nn}
\end{array}
=
\frac {1}{(Pa_0)_{1}}~~
\begin{array}{|c c c|}
(PB)_{11} & \cdots & (PB)_{1n}\\
(PB)_{21} & \cdots & (PB)_{2n}\\
\vdots && \vdots \\
(PB)_{n1} & \cdots & (PB)_{nn}
\end{array}
=
\frac{\det (PB)}{(Pa_0)_{1}}.
\end{displaymath}
But since $P$ is invertible $(Pa_0)_{1}\neq 0$ and $\det P\neq 0$. Thus $D_i=0$ if and only if $\det B=0$. 
\end{proof}

Note that it follows from the proof of Claim ~\ref{claim1} that $\widetilde R_B$ has support on ~$S_\rho$ if $B$ is contained in the facet with normal direction $\rho$. 
Indeed $\widetilde R_B$ survives precisely in the charts corresponding to cones $\tau$ with $\rho$ as one of its generators and in each such chart it has support on the part covering $S_\rho$.

Now let us fix a set $B$ contained in the facet with normal direction $\rho_i$, so that ~\eqref{dennis} is nonvanishing, and compute the action of $\widetilde R_B$ on the pullback of a test form $\phi=\varphi(z)~dz$ of bidegree $(n,0)$. Here $dz$ is just a shorthand notation for $dz_1\wedge\ldots\wedge dz_n$.
 Let $\{\chi_\tau\}$ be a partition of unity on $\mathcal X$ subordinate the cover $\{\mathcal U_\tau\}$. It is not hard to see that we can choose the partition in such a way that the $\chi_\tau$ are circled, that is they only depend on $|t_1|, \ldots, |t_n|$. Now $\widetilde R_B=\sum_{\tau}\chi_\tau \widetilde R_B$. We will start by computing the contribution from our fixed chart $\mathcal U_0$ where $\widetilde R_B$ is realized by ~\eqref{dennis}. 

Since $R_B$ has support at the origin it does only depend on finitely many derivatives of $\varphi$ and therefore to determine $R_B$ it is enough to consider the case when $\varphi$ is a polynomial. We can write $\varphi$ as a finite Taylor expansion,
\begin{equation*}
\varphi=
\sum_{\alpha, \beta}
\frac{\varphi_{\alpha,\beta}(0)}{\alpha ! \beta !} z^\alpha \bar z^\beta, 
\end{equation*}
where $\alpha=(\alpha_1,\ldots,\alpha_n)$, 
$
\varphi_{\alpha,\beta}=
\frac{\partial^{\alpha_1}}{\partial z_1^{\alpha_1}}\cdots
\frac{\partial^{\alpha_n}}{\partial z_n^{\alpha_n}}
\frac{\partial^{\beta_1}}{\partial \bar z_1^{\beta_1}}\cdots
\frac{\partial^{\beta_n}}{\partial \bar z_n^{\beta_n}}\varphi,
$
 and $\alpha!=\alpha_1!\cdots\alpha_n!$,
with pullback to $\mathcal U_0$ given by
\begin{equation*}
\Pi^* \varphi=
\sum_{\alpha, \beta}\frac{\varphi_{\alpha,\beta}(0)}{\alpha ! \beta !}
t^{P\alpha} \bar t^{P\beta}.
\end{equation*}
A computation similar to the one in the proof of Lemma ~\ref{purple} yields
\begin{equation*}
\Pi^* dz=
\det P ~ t^{(P-I)\1} ~dt.
\end{equation*}
Hence $\chi_\tau \widetilde R_B.\Pi^*\phi$ is equal to
\begin{multline*}
K \int
\dbar 
\bigg [\frac{1}{t_i^{n\rho_i\scalar a_0}}\bigg]
\otimes \bigg [\frac{1}{\prod_{j\neq i} t_j^{n\rho_j\scalar a_0}}
\bigg ]
\wedge
\frac{\bar t_i\bar t^{(C-I)\1}~\widehat{d\bar t_i}\wedge e_B}
{\nu(t)^n}
\wedge\\
\chi_\tau(t)
\sum_{\alpha, \beta}
\frac{\varphi_{\alpha,\beta}(0)}{\alpha ! \beta !}t^{P\alpha} \bar t^{P\beta}
t^{(P-I)\1} dt=
K~\sum_{\alpha,\beta} I_{\alpha,\beta}~\wedge
\frac{\varphi_{\alpha,\beta}(0)}{\alpha ! \beta !}
~e_B,
\end{multline*}
where $K= (n-1)! D_1 \det P$ and
\begin{equation}\label{fredag}
I_{\alpha,\beta}
=
\int
\dbar
\bigg [\frac{1}{t_i^{\rho_i\scalar (n a_0-\alpha-\1)+1}}\bigg]
\otimes
[\mu_{\alpha,\beta}]
\frac{ \chi_\tau(t)~\bar t_i^{\rho_1\scalar\beta}}{\nu(t)^n}~
\wedge \widehat{d\bar {t}_i}\wedge dt,
\end{equation}
and where $\mu_{\alpha,\beta}$ is the Laurent monomial in $t_j$ and  $\bar t_j$ for $j\neq i$:
\begin{equation*}
\mu_{\alpha,\beta}=
\prod_{j\neq i}
t_j^{\rho_j\scalar(\alpha+\1-n a_0)-1}~
\bar t_j^{\rho_j\scalar(\beta +B\1-na_0)-1}.
\end{equation*}
Observe that $\rho_j\scalar(\beta +B\1-na_0)-1\geq 0$ so there are no $\bar t_j$ in the denominator. 
Recalling ~\eqref{salt}, we evaluate the $t_i$-integral.
Since $\nu$ and $\chi_\tau$ depend on $|t_1|, \ldots, |t_n|$ 
it follows that
$\frac{\partial^{\ell}}{\partial t_i^\ell}\frac{\chi_\tau}{\nu}|_{t_i=0}=0$ for $\ell \geq 1$ and thus ~\eqref{fredag} is equal to
\begin{equation}\label{skynda}
2\pi i
\int_{\widehat{t_i}}
\frac{\chi_\tau(t)|_{t_i=0}
[\mu_{\alpha,\beta}]}
{\nu(t)^n|_{t_i=0}}~
\widehat{d\bar t_i}\wedge \widehat{dt_i},
\end{equation}
if 
\begin{equation}\label{fearful}
\rho_i\scalar (n a_0-\alpha-\1)+1=1
\end{equation}
and
\begin{equation}\label{captain}
\rho_i\scalar\beta = 0,
\end{equation}
and zero otherwise. 
Moreover, for symmetry reasons ~\eqref{skynda} vanishes unless
\begin{equation}\label{trip}
\rho_j\scalar(\alpha+\1-na_0)-1=
\rho_j\scalar(\beta + B\1-na_0)-1
\end{equation}
for $j\neq i$.
From the discussion just before Theorem ~\ref{main_thm} we know that the facet containing $B$ is compact, which means that its normal vector has nonzero entries. Thus ~\eqref{captain} implies that $\beta=(0,\ldots,0)$. Using the fact that $\rho_i\scalar a=\rho_i\scalar a_0$ for all $a\in B$ we can rewrite the left hand side of ~\eqref{fearful} as $\rho_i \scalar (B\1-1-\alpha-\1)+1$ and thus summarize the conditions ~\eqref{fearful} and ~\eqref{trip} on $\alpha$ as
\begin{equation}\label{john}
P(\alpha+\1)=PB\1.
\end{equation}
But, since $P$ is invertible there exists exactly one $\alpha$ that fulfills ~\eqref{john}, namely $\alpha=(B-I)\1$, which is precisely $\alpha^B-\1$. 
With these values of $\alpha$ and $\beta$ the Laurent monomial $\mu_{\alpha,\beta}$ is nonsingular and so the integrand of ~\eqref{skynda}, 
\begin{equation}\label{johannes}
2\pi i \int_{\widehat{t_i}}
\frac{\chi_\tau(t)|_{t_i=0}~\prod_{j\neq i}|t_j|^{2(\rho_j\scalar (B\1-na_0)-1)}}
{\nu(t)^n|_{t_i=0}}
\widehat{d\bar t_i}\wedge \widehat{d t_i}
\end{equation}
becomes integrable.

To compute $\widetilde R_B.\widetilde\Pi^*\phi$ we want to add contributions from all charts. 
However, $\mathcal U_0$ covers the support of $\widetilde R_B$ except for a set of measure zero, since $\widetilde R_B$ has support on $S_{\rho_i}$, and moreover all integrands that appear are of the form ~\eqref{johannes} and therefore integrable.
Thus $\widetilde R_B.\widetilde \Pi^*\phi$ is equal to 
\begin{multline*}
\int_{\mathcal X}
\sum_\tau 
\widetilde \Pi^*(\dbar|z^A|^{2\lambda}\wedge u_B)
\wedge\chi_\tau\widetilde\Pi^*\phi\Big|_{\lambda=0}=\\
\int_{\mathcal U_0}
\widetilde\Pi^*(\dbar|z^A|^{2\lambda}\wedge u_B)\wedge\widetilde\Pi^*\phi
\Big |_{\lambda=0}=
C_B~\frac{\varphi_{\alpha^B-\1,0}(0)}{\alpha !\beta !}e_B,
\end{multline*}
where 
\begin{equation*}
C_B=
2\pi i K \int_{\widehat{t_i}}
\frac{\prod_{j\neq i}|t_j|^{2(\rho_j\scalar (B\1-na_0)-1)}}
{(\sum_{a\in A}
\prod_{j\neq i}
|t_j|^{2\rho_j\scalar(a-a_0)})^n}
\widehat{d\bar t_i}\wedge \widehat{d t_i}.
\end{equation*}
Hence $R_B$ is of the form ~\eqref{muffel} and the result follows.

\section{General monomial ideals}\label{relative}
If the zero variety of $z^A$ is of positive dimension the computations of $R^{z^A}$ get more involved. Recall that in general $R^{z^A}=R_p+\ldots +R_\mu$, where $p=\codim Y$, $\mu=\min(m,n)$ and $R_k\in \mathcal D'_{0,k}(\mathbb C^n,\Lambda^k E)$. Parts of the top degree term $R_n$ can be computed by the techniques from the proof of Theorem ~\ref{main_thm}. Our method for dealing with the terms of lower degree, though, is to perform the computations outside certain varieties, where some of the coordinates are zero. This amounts to projecting $A$ and brings us back to the more familiar top degree case in a lower dimension. The price we have to pay is that we miss parts of $\mathbb C^n$. More precisely, we will compute the current $R_k$ outside the $(k+1)$-dimensional variety
\begin{equation*}
V_k:=
\bigcup_{{\mathcal I}, |{\mathcal I}|=k+1}\bigcap_{i\in {\mathcal I}} H_i,
\end{equation*}
where $H_i$ denotes the hyperplane $\{z_i=0\}$.
However, it turns out that $R_k$ will not carry any essential information on such ``small'' varieties. To be precise, we have the following lemma, which can be proven analogously to the proof of Lemma 2.2 in ~\cite{A3}. 
\begin{lma}\label{hammer}
Let $h_1,\ldots,h_s$ be a tuple of holomorphic functions and let $Y_h=\{h_1=\ldots h_s=0\}$. Suppose that  $\codim Y_h\cap Y > k$ .
Then the current $ |h|^{2 \lambda_0} R_k$, where $|h|^2=|h_1|^2+\ldots +|h_s|^2$, has an analytic continuation to $\Re\lambda_0\geq -\epsilon$ and
\begin{equation*}
 |h|^{2\lambda_0} R_k|_{\lambda_0=0}=R_k.
\end{equation*}
\end{lma}
It follows, in particular, that to annihilate $R_k$ it suffices to do it outside $V_k$ (or any variety of codimension $k+1$). Indeed $h R_k=0$ outside $V_{k}$ implies that $h R_k=0$.

Before stating our result, a word of notation:
For $\mathcal I=\{i_1,\ldots, i_k\}\subseteq \{1,\ldots,n\}$, let $T_{\mathcal I}$ be the projection
\begin{eqnarray*}
T_{\mathcal I}: \mathbb Z^n &\to & \mathbb Z^k\\
(a_1,\ldots, a_n) &\mapsto & (a_{i_1},\ldots, a_{i_k}).
\end{eqnarray*}
We way that $T_{\mathcal I}(B)$ is essential if $T_{\mathcal I}(B)$ is contained in a facet of $\Gamma^+(T_{\mathcal I}(A))$ and if $T_{\mathcal I}(B)$ spans $\mathbb R^{|\mathcal I|}$.
\begin{thm}
\label{rel_thm}
Let $z^A$, $A\subseteq \mathbb Z_+^n$, be a tuple of monomials in $\mathbb C^n$, and let
\begin{equation*}
R^{z^A}=\sum_{B\subset A} R_B
\end{equation*}
be the corresponding
Bochner-Martinelli residue current.
Then outside $V_{|B|}$,
\begin{equation*}
R_B=
\sum_{\mathcal I\subset\{1,\ldots,n\}, |\mathcal I|=|B|}
R_{B,\mathcal I},
\end{equation*}
where the current $R_{B,\mathcal I}$ vanishes unless $T_{\mathcal I}(B)$ is essential.
Moreover if $T_{\mathcal I} (B)$ is essential and contained in a compact facet of $\Gamma^+(T_{\mathcal I}(A))$, then
\begin{equation}\label{smuggel}
R_{B,\mathcal I}= 
C_{B, \mathcal I}(\eta) \otimes ~ \bigwedge_{i\in \mathcal I} \dbar \Big [\frac{1}{z_i^{\alpha^B_i}} \Big ]\wedge e_B,
\end{equation}
where $\eta$ denotes the $z_i, i\notin \mathcal I$, and $C_{B,\mathcal I}(\eta)$ is a smooth function not identically equal to zero.
\end{thm}
Several remarks are in order. 
First, an immediate consequence is that 
\[\ann R_{B,\mathcal I}=(z_i^{\alpha_i^B})_{i\in \mathcal I}\]
if $R_{B,\mathcal I}$ is of the form ~\eqref{smuggel}, since annihilating such a current clearly is equivalent to annihilating the $\bigwedge_{i\in \mathcal I} \dbar \Big [\frac{1}{z_i^{\alpha^B_i}}\Big ]$ part.
Moreover the support of ~\eqref{smuggel} is the set $\cap_{i\in \mathcal I}\{z_i=0\}$. Note that all the computable $R_{B,\mathcal I}$ have different supports.

\begin{remark}\label{generatorer}
Observe that adding elements to $A$ that lie in any of the non-compact facets of $\Gamma^+(A)$, not contained in any coordinate plane, gives rise to new essential sets. 
For example we can add redundant generators to $(z^A)$ and thus in general $\ann R^{z^A}$ is not independent of the choice of generators as in the case of a discrete zero variety.  
\end{remark}

\begin{remark}
Theorem ~\ref{main_thm} is just a special case of Theorem ~\ref{rel_thm}. Let us say a word about how to see that the currents of lower degree vanish when $Y$ is the origin. This hypothesis 
means precisely that $A$ intersects all axes, which in turn implies that the image of $A$ under any projection $T_{\mathcal I}$, $|\mathcal I|<n$, contains the origin.  However, if $0\in A$, the Newton polyhedron $\Gamma^+(A)$ equals the first orthant and there are no essential sets; note that this corresponds to the case when $f$ contains a nonvanishing function. Thus $R^{z^A}=\sum_{|B|=n}R_B$, where $R_B=R_{B,\{1,\ldots,n\}}$ and Theorem ~\ref{main_thm} follows. Of course, by slightly refined arguments one can see how the currents $R_k, k< \codim Y$ vanish in general.
\end{remark}

\begin{remark}
By Theorem ~\ref{rel_thm} we can extend Theorem ~\ref{main_corollary} to hold for a much larger class of ideals. Recall that a crucial point of the proof of Theorem ~\ref{main_corollary} was the existence of essential sets. If the Newton diagram of $A$ is of dimension $n-1$, though, we can always find essential sets, for example take the vertices of one of the $(n-1)$-dimensional facets, and the proof applies immediately. 
In fact, one can show
that Theorem ~\ref{main_corollary} holds unless $\Gamma(A)$ is not parallel to any of the coordinate planes. Yet, there are ideals for which Theorem ~\ref{rel_thm} does not give enough information to decide whether the inclusion ~\eqref{orm} is strict or not, as we will see in Example ~\ref{notenough}. Still, in this particular case, one can show by explicit computations that the annihilator ideal is strictly included in the ideal and we believe that Theorem ~\ref{main_corollary} holds for monomial ideals in general, although we do not know enough to prove it.
\end{remark}

Let us illustrate Theorem ~\ref{rel_thm} with some simple examples.
\begin{ex}\label{relativeex}
Let
$ A=\{a^1=(6,1), a^2=(3,2), a^3=(2,4)\}$.
There are two essential subsets of $A$, $\{a^1,a^2\}$ and $\{a^2,a^3\}$, with $\alpha^{12}=(9,3)$ and $\alpha^{23}=(5,6)$, respectively. Moreover, $\Gamma^+(T_{\{1\}}(A))$ is the interval $[2,\infty)$ and consequently $\Gamma(T_{\{1\}}(A))=\{2\}$. Thus the only set such that its image under $T_{\{1\}}$ is essential is $\{a^3\}$, with $\alpha^3=a^3$, and according to Theorem ~\ref{rel_thm} $\ann R_{\{a^3\},\{1\}}=(z_1^2)$. Similarly, projecting $A$ on the second axis yields one current, $R_{\{a^1\},\{2\}}$, with annihilator $(z_2)$. Altogether we get 
\[\ann R^{z^A}=(z_1^9,z_2^3)\cap (z_1^5,z_2^6)\cap (z_1^2)\cap (z_2),\]
that is equal to
$(z_1^9z_2,z_1^5z_2^3,z_1^2z_2^6)$, 
see Figure 7.
\begin{figure}\label{relex}
\begin{center}
\includegraphics{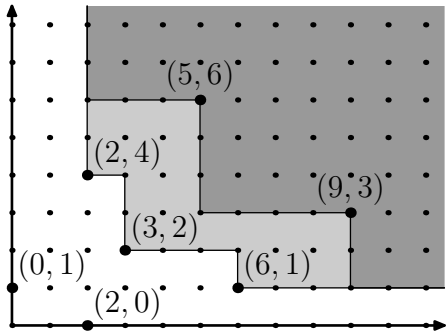}
\caption{The ideals $\ann R^{z^A}$ (dark grey) and $(z^A)$ (light grey) in Example ~\ref{relativeex}}
\end{center}
\end{figure}
Observe, apropos of Remark ~\ref{generatorer}, that adding a point to $A$ in any of the noncompact facets gives a new essential set and thereby essentially changes $R^{z^A}$.
\end{ex}
In view of Example ~\ref{relativeex} it should be clear that Theorem ~\ref{rel_thm} actually gives a complete description of $\ann R^{z^A}$ in case $n=2$, provided we choose a minimal set of generators (or at least avoid to pick redundant generators from the unbounded facets of $\Gamma^+(A)$). 

\begin{ex}\label{notenough}
Let $I$ be the ideal $(z^A)$, where 
$
A=\{a^1=(1,0,1), a^2=(0,1,1)\} \subset \mathbb Z^3.
$
The codimension of $\{z^a=0\}$ is $1$ and thus $I$ is not a complete intersection (nor can be defined by one). Note that the set $A$ is to small to be essential, whereas the image of $A$ under any projection to $\mathbb Z^2$ is, as shown in Figure 8. Still, Theorem ~\ref{rel_thm} gives the annihilator ideal only for one of the corresponding currents, namely $\ann R_{A,\{1,2\}}=(z_1, z_2)$. In both of the other cases the projection of $A$ lies in a noncompact facet of the Newton polyhedron. Furthermore, projecting $A$ to $\mathbb Z$ yields the currents $R_{\{a^1\},\{3\}}$ and $R_{\{a^2\},\{3\}}$, both with annihilator $(z_3)$. Observe that the intersection of the computable currents is precisely $I$. Thus we have found an example of an non-complete intersection where Theorem ~\ref{rel_thm} does not give enough information to decide whether the inclusion ~\eqref{orm} is strict or not. In this simple example, however, it is easy to compute the remaining parts of $R^{z^A}$ and see that the inclusion is indeed strict.
\begin{figure}\label{noten}
\begin{center}
\includegraphics{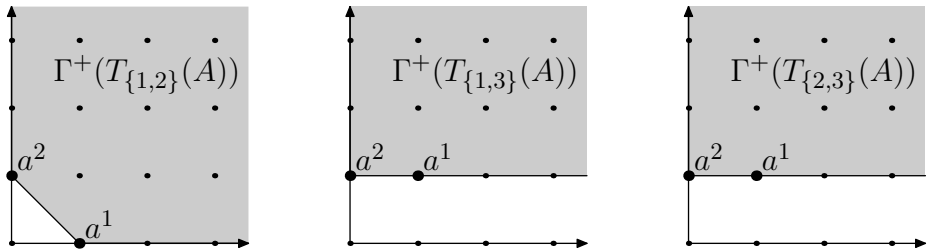}
\caption{The image of $A$ under the various projections to $\mathbb Z^2$ in Example ~\ref{notenough}}
\end{center}
\end{figure}
\end{ex}

\begin{proof}[Proof of Theorem ~\ref{rel_thm}]
We start by considering the term of top degree,
\begin{equation*}
R_n=\sum_{B\subset A, |B|=n} R_B,
\end{equation*}
for which the result follows easily from the proof of Theorem ~\ref{main_thm}. To see this, observe first that the proof of Claim ~\ref{claim1}
does not depend on the codimension of $Y$. Thus we conclude that $R_B=0$ unless $B$ is essential.

Next, suppose that $B$ is contained in a compact facet $F_B$ of $\Gamma^+$ with normal direction $\rho_i$. 
As in the proof of Theorem ~\ref{main_thm} let $\mathcal U_0$ be a chart parametrized by $\Pi$, determined by the cone $\tau_0$ that has $\rho_i$ as its $i$th generator. 
Recall from the proof that the support of $\widetilde R_B$ in $\mathcal U_0$ is given by $\{t_i=0\}$. That $F_B$ is compact means precisely that all entries of $\rho_i$ are strictly positive, which implies that $\Pi(\{\rho_i=0\})=\{0\}$. Consequently, when computing $\widetilde R_B$ in $\mathcal U_0$, we only need to consider it acting on test forms $\phi=\varphi~dz$, where $\varphi$ is a polynomial. Hence the rest of of the proof of Theorem ~\ref{main_thm} applies, and we get that $R_B=R_{B,\{1,\ldots,n\}}$ is of the form ~\eqref{muffel} that is equivalent to ~\eqref{smuggel} in case $k=n$.

We will compute the terms of lower degree by looking outside certain coordinate planes, which will correspond to projections of $A$. More precisely, to determine $R_k$ we will look where $n-k$ of the $z_i$ are nonzero. To do this let us fix ${\mathcal I}=\{i_1,\ldots,i_{k}\}\subseteq \{1,\ldots,n\}$ and let $M_{\mathcal I}$ be the set where $z_i$ is nonvanishing if $i\notin {\mathcal I}$, that is
\[
M_{\mathcal I}= (\bigcup_{i\notin {\mathcal I}} H_i)^C.
\]
Denote the $z_i, i\in {\mathcal I}$, by $\zeta$ and the $z_i, i\notin {\mathcal I}$, by $\eta$ and write $z^a=\zeta^{a_\zeta}\eta^{a_\eta}$, where $a_\zeta$ and $a_\eta$ are the images of $a$ under $T_{\mathcal I}$ and $T_{\mathcal I^C}$, respectively. Let $A_\zeta$ and $A_\eta$ denote the corresponding images of $A$, and let $\phi$ be a test form of bidegree $(n,k)$ with (compact) support in $M_{\mathcal I}$. Now $R_k$ acting on $\phi$ is the analytic continuation to $\lambda=0$ of
\begin{equation*}
\int
\dbar|z^A|^{2\lambda}\wedge 
\frac{s\wedge(\dbar s)^{k-1}}{|z^A|^{2k}}\wedge \phi(z),
\end{equation*}
that is equal to a sum, taken over $B$ such that $|B|=k$, of terms 
\begin{equation}\label{luke}
\int_{\eta}\int_{\zeta}
\dbar|\zeta^{A_\zeta}\eta^{A_\eta}|^{2\lambda}\wedge 
\frac{s_B\wedge(\dbar s_B)^{k-1}}{|z^A|^{2k}}
\wedge \varphi(\zeta,\eta)~d\zeta\wedge d\bar\eta\wedge d\eta.
\end{equation}
It is easily checked that $R_k$ vanishes unless $\phi$ is of the form $\varphi(\zeta, \eta) ~d\bar\eta\wedge d\eta\wedge d\zeta$. We can now compute the inner integral of ~\eqref{luke} as in the top degree case (with $A^\zeta$ in $\mathbb C^k_\zeta$). Indeed, since $\eta$ is nonvanishing, we can regard $z^A$ as the monomials $\zeta^{A_\zeta}$ times the parameters $\eta^{a_\eta}$. It follows that, at $\lambda=0$, ~\eqref{luke} vanishes unless $T_{\mathcal I}(B)$ is essential, and moreover, if $T(B)$ is contained in a compact facet of $\Gamma(T(A))$, then the inner integral is equal to
\begin{equation*}
C_{B,\mathcal I}(\eta) \otimes ~\dbar \bigg [\frac{1}{\zeta_1^{\alpha^B_1}} \bigg ]\wedge\ldots\wedge
 \dbar \bigg[\frac{1}{\zeta_k^{\alpha^B_k}}\bigg]\wedge e_B\wedge \varphi(\zeta,\eta)~d\zeta,
\end{equation*}
where $C_{B,\mathcal I}$ depends smoothly on $\eta$.

In other words, if we let $R_{B,\mathcal I}$ be defined by ~\eqref{luke}, meaning that its action on a test form $\phi$ is the value of ~\eqref{luke} at $\lambda=0$), then $R_{B,\mathcal I}$ is of the form ~\eqref{smuggel}.

When looking in $M_{\mathcal J}$ for each index set ${\mathcal J}$ of cardinality $k$ we miss
\begin{multline*}
(\bigcup_{{\mathcal J}, |{\mathcal J}|=k}M_{\mathcal J})^C=
(\bigcup_{{\mathcal J}, |{\mathcal J}|=k}(\bigcup_{i\notin {\mathcal J}} H_i)^C)^C=\\
\bigcap_{{\mathcal J}, |{\mathcal J}|=k}\bigcup_{i\notin {\mathcal J}} H_i=
\bigcup_{{\mathcal J}, |{\mathcal J}|=k+1}\bigcap_{i\in {\mathcal J}} H_i, 
\end{multline*}
that is precisely $V_k$. Clearly each current $R_{B,\mathcal I}$ extends to $\bigcup_{{\mathcal J}, |{\mathcal J}|=k}M_{\mathcal J}$. In fact $R_{B,\mathcal I}$ has support only in $M_{\mathcal I}$. Thus outside $V_k$ we have $R_k=\sum R_{B,\mathcal I}$, where the $R_{B,\mathcal I}$ are of the desired form and we are done.
\end{proof}

{\bf Acknowledgement:} The author would like to thank Mats Andersson for interesting discussions on the topic of this paper and for valuable comments on preliminary versions.

\def\listing#1#2#3{{\sc #1}:\ {\it #2},\ #3.}


\begin{thebibliography}{9999}

\bibitem{A}\listing{M.\ Andersson}
{Residue currents and ideals of holomorphic functions.}
{Bull.\ Sci.\ Math. {\bf 128} (2004) no. 6 481--512}

\bibitem{A2}\listing{M.\ Andersson}
{Integral representations with weights I}
{Math.\ Ann.\ {\bf 326} (2003),  no. 1, 1--18}

\bibitem{A3}\listing{M.\ Andersson}
{The membership problem for polynomial ideals in terms of residue currents}
{Ann.\ Inst.\ Fourier (to appear)}

\bibitem{BGVY}\listing{C.\ A.\ Berenstein \& R.\ Gay \& A.\ Vidras \&
A.\ Yger} {Residue currents and Bezout identities}
{Progress in Mathematics
{\bf 114} Birkh\"auser Yerlag (1993)}

\bibitem{BS}\listing{J.\ Brian\c{c}on, H.\ Skoda }
{Sur la cl\^oture int\'egrale d'un id\'eal de germes de fonctions holomorphes en un point de $\mathbb C^n$}
{C.\ R.\ Acad.\ Sci.\ Paris S\'er.\ A {\bf 278} (1974) 949--951}

\bibitem{DS}\listing{A.\ Dickenstein  \& C.\ Sessa}
{Canonical representatives in moderate cohomology}
{Invent. Math. {\bf 80} (1985),  417--434}

\bibitem{KKMS}\listing{G.\ Kempf \& F.\ Knudsen \& D.\ Mumford \& B.\ Saint-Donat}
{Toroidal Embeddings I} 
{Lecture Notes in Mathematics {\bf 339} Springer Verlag, New York, 1973}

\bibitem{K}\listing{A.\ G.\ Khovanskii}
{Newton polyhedra and torodial varieties}
{Funct.\ Anal.\ Appl. {\bf 11} (1978), 289 -- 295}

\bibitem{MS}\listing{E.\ Miller \& B.\ Sturmfels}
{Combinatorial commutative algebra} 
{Graduate Texts in Mathematics {\bf 227} Springer-Verlag, New York, 2005}

\bibitem{P}\listing{M.\ Passare}
{Residues, currents, and their relation to ideals of holomorphic functions} 
{Math.\ Scand.\ {\bf 62} (1988), no. 1, 75--152}

\bibitem{PTY}\listing{M.\ Passare \& A.\ Tsikh \&  A.\ Yger}
{Residue currents of the Bochner-Martinelli type}
{Publ.\ Mat.  {\bf 44} (2000), 85--117}

\bibitem{T}\listing{B.\ Teissier}{Vari\'et\'es polaires. II. Multiplicit\'es polaires, sections planes, et conditions de Whitney \emph{Algebraic geometry (La R\'abida, 1981)}}
{Lecture Notes in Mathematics {\bf 961} Springer Verlag, Berlin, 1982, pp 314--491}

\bibitem{V}\listing{A.\ N.\ Varchenko}{Newton polyhedra and estimating of oscillating integrals}{Funct.\ Anal.\ Appl. {\bf 10} (1976), 175 -- 196}

\bibitem{Z}\listing{G.\ Ziegler}
{Lectures on polytopes}
{Graduate Texts in Mathematics {\bf 152} Springer Verlag, New York, 1995}

\end{thebibliography}
\end{document}